\documentclass[12pt]{article}
\usepackage{amsmath}
\usepackage{amssymb}
\begin{document}
\newtheorem{theorem}{Theorem}[section]
\newtheorem{remark}[theorem]{Remark}
\newtheorem{mtheorem}[theorem]{Main Theorem}
\newtheorem{observation}[theorem]{Observation}
\newtheorem{proposition}[theorem]{Proposition}
\newtheorem{lemma}[theorem]{Lemma}
\newtheorem{mlemma}[theorem]{Main Lemma}
\newtheorem{note}[theorem]{}
\newtheorem{1killlemma}[theorem]{First Killing-Lemma}
\newtheorem{2killlemma}[theorem]{Second Killing-Lemma}
\newtheorem{corollary}[theorem]{Corollary}
\newtheorem{notation}[theorem]{Notation}
\newtheorem{example}[theorem]{Example}
\newtheorem{definition}[theorem]{Definition}

\renewcommand{\labelenumi}{(\roman{enumi})}
\newcommand{\dach}[1]{\hat{\vphantom{#1}}}
\numberwithin{equation}{section}
\def\Z{{ \mathbb Z}}
\def\E{{ \mathbb E}}
\def\H{{ \mathbb H}}
\def\N{{ \mathbb N}}
\def\DD{ \mathbb D}
\def\DDD{ \mathfrak D}
\def\GG{ \mathbb G}
\def\BBB{ B_\BB}
\def\D{{\hat{D}}}
\def\Q{{\mathbb Q}}
\def\G{\hat{G}}
\def\C{\hat{C}}
\def\T{{\cal T}}
\def\FF{{\mathfrak  F}}
\def\PP{{\mathfrak  P}}
\def\FFF{{\mathfrak  F}^*}
\def\C{{\mathfrak C}}
\def\X{{\mathfrak X}}
\def\Y{{\mathfrak Y}}
\def\RX{R\langle X \rangle}
\def\GX{\langle X \rangle}
\def\bm{\bar{m}^p}
\def\bmq{\bar{m}^q}
\def\bma{\langle \bar{m}^p_{\bar{a}}: \bar{a} \in u^p\rangle}
\def\ba{\bar{a}^p}
\def\baq{\bar{a}^q}
\def\bap{\langle a^p_l : l < l^p \rangle}
\def\K{{\mathfrak K}}
\def\R{\widehat{R}}
\def\wZ{\widehat{\Z}}
\def\F{\widehat{F}}
\def\G{\widehat{G}}
\def\T{{\cal T}}
\def\B{\widehat{B}}
\def\BC{\widehat{B_C}}
\def\BCC{\widehat{B_\C}}
\def\restr{\restriction}
\def\Aut{{\rm Aut\,}}
\def\Map{{\rm Map\,}}
\def\Im{{\rm Im\,}}
\def\ker{{\rm ker\,}}
\def\lg{{\rm lg\,}}
\def\br{{\rm br\,}}
\def\inf{{\rm inf\,}}
\def\sup{{\rm sup\,}}
\def\Br{{\rm Br\,}}
\def\Yphi{Y_{[\phi]}}
\def\Ypsi{Y_{[\psi]}}
\def\Xphi{X_{\tilde{\phi}}}
\def\Xpsi{X_{\tilde{\psi}}}
\def\a{\alpha}
\def\abar{\overline{\alpha}}
\def\aa{{\bf a}}
\def\to{\rightarrow}
\def\arr{\longrightarrow}
\def\sigmaa{{\bf \Sigma_a}}
\def\End{{\rm End\,}}
\def\Ult{{\rm Ult\,}}
\def\Ines{{\rm Ines\,}}
\def\Hom{{\rm Hom\,}}
\def\Fin{{\rm Fin\,}}
\def\restr{\upharpoonright}
\def\Ext{{\rm Ext}\,}
\def\Hom{{\rm Hom}\,}
\def\End{{\rm End}\,}
\def\Aut{{\rm Aut}\,}
\def\ker{{\rm ker}\,}
\def\Ann{{\rm Ann}\,}
\def\defe{{\rm def}\,}
\def\rk{{\rm rk}\,}
\def\crk{{\rm crk}\,}
\def\nuc{{\rm nuc}\,}
\def\Dom{{\rm Dom}\,}
\def\Im{{\rm Im}\,}
\def\Yphi{Y_{[\phi]}}
\def\Ypsi{Y_{[\psi]}}
\def\Xphi{X_{\tilde{\phi}}}
\def\Xpsi{X_{\tilde{\psi}}}
\def\abar{\overline{\alpha}}
\def\ra{\rightarrow}
\def\arr{\longrightarrow}
\def\iff{\Longleftrightarrow}
\def\mm{{\mathfrak m}}
\def\X{{\mathfrak X}}
\def\Diam{\diamondsuit}

\title{{\sc Groups isomorphic to all their non--trivial normal subgroups}
\footnotetext{This work is supported by Project
No. G-545-173.06/97 of the German-Israeli
Foundation for Scientific Research \& Development\\
AMS subject classification:\\
primary: 20E06, 20F14, 20F22  \\
secondary: 03C65, 03C68, 03C98 \\
Key words and phrases: $\kappa$--class, homogeneous and universal structure,
HNN extension, free product \\ with amalgamation\\
{[GPSh:740]}  in Shelah's list of publications.
}}

\date{}
\author{ R\"udiger G\"obel, Agnes T. Paras and Saharon Shelah}
\maketitle
\begin{abstract}
In answer to a question of P. Hall, we supply another construction of a
group which is isomorphic to each of its non--trivial normal
subgroups.
\end{abstract}

\section{Introduction}
In the early seventies, the following question was asked by Philip Hall
of J.C. Lennox (see \cite{L}) and appeared later in \cite{KN} in the
form :
\begin{quote}
Must a non--trivial group, which is isomorphic to each of its non--trivial
normal subgroups, be either free of infinite rank, simple or infinite
cyclic?
\end{quote}
This was answered in the negative by Obraztsov in \cite{O}, where he used
the technique of graded diagrams developed by Ol'shanskii in \cite{Ol}.

In this note we provide a different construction of a non--trivial group
which is isomorphic to each of its non--trivial normal subgroups. This
construction makes special use of the notion of homogeneous and universal
structures in model theory and well--known facts about HNN extensions and free
products with amalgamation. In our view, the following construction yields a
transparent proof of the existence of such groups and shows how basic ideas
from model theory can be used to resolve questions in other areas of
mathematics. While Obraztsov's proof is based on the deep theory of the
geometry of defining relations in groups, our approach is accessible to a
graduate student with basic knowledge of model theory and group theory. Our
main result is the following:
\begin{quote}
{\bf Main Theorem.} Let $\kappa$ be a cardinal such that
$\kappa = \kappa ^{< \kappa}$. Then there exists a group $G$ of cardinal
$\kappa$ with a descending principal series
$$G \supset G_1 \supset \cdots \supset G_n \supset \cdots$$
such that $\bigcap _{n < \omega} G_n = 1$ and any normal subgroup
of $G$ is some $G_i$, which is isomorphic to $G$.
\end{quote}

\section{Homogeneous and universal structures}
In this section we consider a class of models $\cal M$ and give a brief
summary of its relevant properties. Notation and terminology observed here
follow that in \cite{B}.

Let $\cal M$ be a collection of relational structures in the language $L$ of
groups with constant symbol $1$ (whose value in any model will be the group
identity), the axioms of groups and an additional sequence
$P_n$ ($n \in \omega$) of unary predicates having the following properties:
If $A$ is a group and $P_nA$ denotes $P_n$ applied to $A$, then
$P_nA$ is a normal subgroup of $A$ (denoted $P_nA \triangleleft A$)
and $P_{n+1}A \subseteq P_nA$. In terms of first order
language, this is to say that
\begin{enumerate}
\item[(1)] There exists $1 \in P_nA$.
\item[(2)] If $x, y \in P_nA$, then $xy^{-1} \in P_nA$.
\item[(3)] If $x \in P_nA$ and $z \in A$, then $z^{-1}xz \in P_nA$.
\item[(4)] If $x \in P_{n+1}A$, then $x \in P_nA$.
\end{enumerate}
Thus the class $\cal M$ of all such $L$-models is
\[ {\cal M} = \{ (P_{n}A)_{n < \omega} : P_{n}A \mbox{ is a group },
P_{n}A \triangleleft P_{0}A \mbox{ and } P_{n+1}A \subseteq P_{n}A \,\,
(n < \omega) \} \]
To simplify notation, we let $P_0A = A$ and abbreviate the sequence
$(P_nA)_{n < \omega}$ of groups as $(P_nA)$.

Let $(P_nA)$, $(P_nB) \in \cal M$. The relational structure $(P_nA)$ is said
to be a {\em substructure} of $(P_nB)$ (equivalently, $(P_nB)$ is an
{\em extension} of $(P_nA)$) if $P_nA = P_nB \cap P_0A$ for all $n < \omega$.
We denote this by $(P_nA) \leq (P_nB)$. Let
$\varphi : P_0A \rightarrow P_0B$ be a map. We say that $\varphi$ is a
{\em homomorphism}
of $(P_nA)$ into $(P_nB)$ if $\varphi$ is a group homomorphism and
$\varphi (P_nA) = P_nB \cap \varphi (P_0A)$ ($n < \omega$). A homomorphism
$\varphi$ is called an
{\em embedding} of $(P_nA)$ into $(P_nB)$ if $\varphi$ is one-to-one. Hence if
$(P_nA) \leq (P_nB)$, then the identity map
$id : (P_nA) \rightarrow (P_nB)$ is an embedding. If $\varphi$
is one-to-one and onto, then $\varphi$ is called an {\em isomorphism} and
$(P_nA)$ is said to be {\em isomorphic} to $(P_nB)$. This we denote by
$(P_nA) \cong (P_nB)$. We refer to the cardinal of $P_0A$ as the
{\em cardinal of} $(P_nA)$.

\begin{proposition}\label{kclass}
The following properties hold in $\cal M$ for an infinite cardinal $\kappa$:
\begin{enumerate}
\item[(I)] $\cal M$ contains structures of arbitrarily large cardinality.
\item[(II)] If $(P_nA) \in \cal M$ and $(P_nA) \cong (P_nB)$, then
$(P_nB) \in \cal M$.
\item[(III)] If $(P_nA_i)$, $(P_nB) \in \cal M$, there exist
$(P_nC) \in \cal M$ and embeddings
$$\varphi _i : (P_nA_i) \hookrightarrow (P_nC) \,\,(i = 1,2).$$
\item[(IV)] (Amalgamation Property) Let
$(P_nA)$, $(P_nB_1)$, $(P_nB_2) \in \cal M$ such that
$$\varphi _i : (P_nA) \hookrightarrow (P_nB_i) \,\, (i = 1,2)$$ are embeddings.
Then there is some $(P_nC) \in \cal M$ and embeddings
$$\gamma _i : (P_nB_i) \hookrightarrow (P_nC) \,\, (i = 1,2) \,\,such \,\,
that \,\, \gamma _1 \circ \varphi _1 = \gamma _2 \circ \varphi _2 \,.$$
\item[(V)] The union of any chain of structures
$(P_nA) \in \cal M$ is still in $\cal M$.
\item[(VI$_{\kappa}$)] If $(P_nA) \in \cal M$ and
$(P_nC) \leq (P_nA)$ such that $(P_nC)$ is of cardinal $< \kappa$, then
there exists $(P_nB) \in \cal M$ of cardinal $< \kappa$ such that
$(P_nC) \leq (P_nB) \leq (P_nA)$.
\item[(VII$_{\kappa}$)] Each $(P_nA) \in \cal M$ has $< \kappa$
relations.
\end{enumerate}
\end{proposition}

\noindent
{\bf Proof.} (II) This is clear by the definition of a homomorphism between
elements of $\cal M$. (V) is obvious. (VII$_{\kappa}$) is clear, since, by
definition, each element of $\cal M$ has $\aleph_0$ relations.

(I) Let $P_0A = F_{\alpha}$ be a free group with rank $\alpha$ and
$P_nA = F_{\alpha}^{(n)}$ be the $n$th derived group of $F_{\alpha}$.
Then $(P_nA) \in \cal M$ with arbitrarily large cardinality $\alpha$.

(III) Take, for instance, $P_nC = P_nA_1 \oplus P_nA_2$ ($n < \omega$) and
the natural injections $\varphi _i : P_0A_i \rightarrow P_0C$ ($i = 1,2$).

(IV) Consider the free product with amalgamation
$P_0C = P_0B_1 \ast_{P_0A} P_0B_2$ and the identity map
$id_i : P_0B_i \rightarrow P_0C$. Then
$id_1 \circ \varphi _1 = id_2 \circ \varphi _2$,
since $\varphi _1(a) = \varphi _2(a)$ ($a \in P_0A$)
when viewed as elements of $P_0C$.
Define $P_nC = \langle P_nB_1, P_nB_2 \rangle^{P_0C}$, the normal subgroup
generated by $\langle P_nB_1, P_nB_2 \rangle$ in $P_0C$. Clearly
$(P_nC) \in \cal M$. It suffices to verify that
$P_nC \cap P_0B_i = P_nB_i$ ($n < \omega$). Let
$G_n^* = P_0B_1/P_nB_1 \ast_{P_0A/P_nA} P_0B_2/P_nB_2$ be a free product with
amalgamation and
$\bar{\varphi _{ni}}: P_0A/P_nA \hookrightarrow P_0B_i/P_nB_i$ ($i = 1,2$)
be the homomorphisms defined by
$\bar{\varphi _{ni}}(aP_nA) = \varphi _i(a)P_nB_i$.
Consider the canonical epimorphisms
$\pi _{ni}: P_0B_i \rightarrow P_0B_i/P_nB_i$ ($n < \omega$, $i = 1,2$). Then
the map $\pi _{n1} \cup \pi _{n2} : P_0B_1 \cup P_0B_2 \rightarrow G_n^*$
extends to a homomorphism $\pi : P_0C \rightarrow G_n^*$ such that
$\pi\restriction _{P_0B_i} = \pi _{ni}$. If
$x \in P_0B_i \setminus P_nB_i$, then $\pi _{ni}(x) \neq 1$ and so
$\pi (x) \neq 1$. If $x \in P_nB_1 \cup P_nB_2$, then $\pi (x) = 1$. Thus
ker$\pi = \langle P_nB_1, P_nB_2 \rangle^{P_0C} = P_nC$ and so
$P_nC \cap P_0B_i = P_nB_i$ ($i = 1,2$). Thus $(P_nB_i) \leq (P_nC)$.

(VI$_{\kappa}$) Note that $(P_nC)$ is itself in $\cal M$, hence take
$(P_nB) = (P_nC)$. $\hfill \square$

\bigskip
The preceding properties (I)--(VII)$_{\kappa}$ show that $\cal M$ is a
$\kappa$--class, where $\kappa$ is an infinite cardinal. Note that
if ${\cal M}^+$ is defined to consist of all elements in $\cal M$ such that
the first $n$ terms of the sequence are all equal (for a fixed $n$),
then ${\cal M}^+$ is still a $\kappa$--class.

A relational structure $(P_nA) \in \cal M$ is said to be
($\cal M, \kappa$)--{\em homogeneous} if given any two
$(P_nB_1)$, $(P_nB_2) \in \cal M$ of cardinal $< \kappa$ which are
substructures of $(P_nA)$ and an isomorphism
$\varphi : (P_nB_1) \rightarrow (P_nB_2)$, then there is an isomorphism
$\gamma : (P_nA) \rightarrow (P_nA)$ such that
$\gamma\restriction_{(P_nB_1)} = \varphi$. A structure $(P_nA)$ is said to
be $\cal M$--{\em homogeneous} if $(P_nA)$ is of cardinal $\kappa$ and
($\cal M, \kappa$)--homogeneous. $(P_nA)$ is said to be
$\cal M$--{\em universal}, if $(P_nA)$ is of cardinal $\kappa$ and
given any $(P_nB) \in \cal M$ of cardinal $< \kappa ^+$, there is an
embedding
of $(P_nB)$ into $(P_nA)$. By J\'onsson's Theorem (see \cite{B}, p.213 or
\cite{DG}), if $\kappa$ is chosen such that $\kappa = \kappa ^{< \kappa}$
(e.g., $\kappa$ is regular and either $\kappa$ is a limit beth number or
the GCH holds), then $\cal M$ contains an $\cal M$--homogeneous,
$\cal M$--universal structure of cardinal $\kappa$ which is unique up to
isomorphism. We denote this structure by $(P_nG)$.

\section{Normal subgroups of $P_0G$}

We state and verify here special properties of $P_0G$ and use them to obtain
the required group.
\begin{proposition}\label{isomorphic}
Let $(P_nG)$ be $\cal M$--homogeneous, $\cal M$--universal of cardinal
$\kappa$. Then $P_nG \cong P_0G$ ($n < \omega$).
\end{proposition}

\noindent
{\bf Proof.} Let $M = (P_nG)$. Define
$M^n = (P_nG, P_nG, \cdots , P_nG, P_{n+1}G, P_{n+2}G, \cdots )$,
where the first $n+1$ terms in the sequence are all equal to $P_nG$ and the
succeeding terms agree with the corresponding terms in $M$, i.e., if
$k > 1$, the $(n+k)$th term of $M^n$ is the $(n+k)$th term of $M$.
Define $${\cal M}^n = \{ (P_nA) \in {\cal M} : \mbox{ the first $n+1$ terms
of $(P_nA)$ are equal }\}$$ which is a $\kappa$--class containing $M^n$.
Since $M$ is $\cal M$--homogeneous, $\cal M$--universal and $M^n \leq M$,
it follows that $M^n$ is ${\cal M}^n$--homogeneous and ${\cal M}^n$--universal.

Consider $M^- = (P_nG, P_{n+1}G, \cdots )$. We show that $M^-$ is
$\cal M$--homogeneous and $\cal M$--universal. For then, by
J\'onsson's Theorem, $M$ must be isomorphic to $M^-$ and so
we have shown that $P_0G \cong P_nG$.

Suppose $X = (P_nA) \in \cal M$ of cardinal $< \kappa ^+$. Define
$$X^0 = (P_0A, P_0A, \cdots , P_0A, P_0, P_1A, P_2A, \cdots )$$ to be an
element
of ${\cal M}^n$, where the first $n$ terms are all equal to $P_0A$ and
the succeeding terms of the chain $X^0$ are $P_nA$ ($n < \omega$).
Since $M^n$ is ${\cal M}^n$--universal,
there exists an embedding $\varphi : X^0 \hookrightarrow M^n$. Hence
$\varphi \restriction _{X} : X \hookrightarrow M^-$ is also an embedding, and so
$M^-$  is $\cal M$--universal.

Let $X = (P_nA)$ and $Y = (P_nB)$ be substructures of $M^-$ and
$\varphi : (P_nA) \rightarrow (P_nB)$ be an isomorphism. We show that there
exists an automorphism $\gamma$ of $M^-$ such that
$\gamma\restriction _X = \varphi$. As in the
preceding paragraph, define the
elements $X^0$ and $Y^0$, which are substructures of $M^n$.
Then $\varphi: X^0 \rightarrow Y^0$ is also an isomorphism. Since
$M^n$ is ${\cal M}^n$--homogeneous, there exists an automorphism $\rho$ of
$M^n$ such that $\rho\restriction _{X^0} = \varphi$. It is clear that
$\gamma = \rho\restriction _{M^-}$ is the required automorphism of
$M^-$ which extends
$\varphi$. Thus we have shown that $M^-$ is $\cal M$--homogeneous.
$\hfill \square$

\bigskip
The following results will establish that if $x \in P_nG \setminus P_{n+1}G$,
then $\langle x \rangle ^{P_0G} = P_nG$. Thus if $N$ is a normal subgroup
of $P_0G$ which is not contained in $\bigcap _{n < \omega} P_nG$, then
$N = P_kG$ for some $k$.

\begin{lemma}\label{intersection}
Let $H$ be a normal subgroup of a group $G$ and $x, y \in G$
with $|x| = |y|$. Assume that for all natural numbers $k$,
$x^k \in H$ if and only if $y^k \in H$. If $G^* = \langle G, t \rangle$
is the HNN extension of $G$ such that $x^t = y$, then
$H^{G^*} \cap G = H$.
\end{lemma}

\noindent
{\bf Proof.} Note that $H^{G^*} = H^{\langle t \rangle}$, since $H$ is
normal in $G$. Clearly, $H \subseteq H^{\langle t \rangle}$. Let
$g \in H^{\langle t \rangle} \cap G$. Then
$g = h_1^{t^{n_1}}\cdots h_k^{t^{n_k}}$, for some $h_i \in H$ and
integers $n_i$. Without loss of generality, assume that
$h_i \neq 1$ and $g$ cannot be written as a product of fewer than $k$ such
conjugates $h_i^{t^{n_i}}$. If $h_i =  x^k$, for some $k$ and
$n_i \geq 0$, then replace $h_i^{t^{n_i}}$ by $y^k$.  Similarly, if
$h_i = y^k$ and $n_i \leq 0$, replace $h_i^{t^{n_i}}$ by $x^k$. Thus we can
further assume that, if $h_i^{t^{n_i}} \in H$, then $n_i = 0$.
If all the $n_i = 0$, then it is clear that $g \in H$. Now suppose,
without loss of generality, that $n_1 \neq 0$. So when $g$ is written
in normal form, it will have length $\geq 1$
(length of an element $g$ refers to the number of occurrences of $t^{\pm 1}$ in its normal form).
But the elements of $G$ have length $0$ (see \cite{LS}, p. 178). Thus, all
the $n_i$ must be $0$. $\hfill \square$

\begin{theorem}\label{conjugates}
Suppose $x, y \in P_nG$ such that $|x| =|y|$ and
$(\langle x \rangle \cup \langle y \rangle ) \cap P_{n+1}G = 1$.
Then there exists $t \in P_0G$ such that $x^t = y$.
\end{theorem}

\noindent
{\bf Proof.} Consider $P_0S = \langle x, y \rangle$ and
$P_nS = P_nG \cap P_0S$ ($n < \omega$). Let $P_0T = \langle x, y, t \rangle$
be the HNN extension of $P_0S$ and $P_nT = P_nS^{P_0T}$ be the normal
subgroup generated by $P_nS$ in $P_0T$.
Then $(P_nS), (P_nT) \in \cal M$ and $(P_nS) \leq (P_nG)$. By
Lemma \ref{intersection}, $(P_nS) \leq (P_nT)$.
Since $(P_nG)$ is $\cal M$--universal, there exists an embedding
$\varphi$ from  $(P_nT)$ to $(P_nG)$. Let
$P_nU = \varphi (P_nT)$ and $P_nV = \varphi (P_nS)$.
Then $(P_nV) \leq (P_nU) \leq (P_nG)$ and $(P_nS) \cong (P_nV)$ via $\varphi$.
Since $(P_nG)$ is $\cal M$--homogeneous, there exists an isomorphism
$\gamma : (P_nG) \rightarrow (P_nG)$ such that
$\gamma\restriction _{(P_nS)} = \varphi$.
Now $\varphi (t) \in P_0U$ and
$\varphi (x^t) = \varphi (x)^{\varphi (t)} = \varphi (y)$
imply that there exists
$t'' \in P_0G$ such that $\gamma (t'') = \varphi (t)$. Since $x, y \in P_0S$,
$\gamma (x^{t''}) = \gamma (x)^{\gamma (t'')} = \varphi (x)^{\varphi (t)} =
\varphi (y) = \gamma (y)$. Since $\gamma$ is an
isomorphism, $x^{t''} = y$ for some $t'' \in P_0G$. Hence $x$ and $y$ are
conjugate in $P_0G$ via the element $t''$. $\hfill \square$

\bigskip
Recall that in Proposition \ref{kclass} (I), we defined $(P_nA)$ to be
$P_nA = F_{\alpha}^{(n)}$, where $F_{\alpha}$ is a free group of infinite rank
$\alpha$. Since $P_nA/P_{n+1}A$ is free abelian of rank $\alpha$, it is clear
that there exist infinitely many elements $x \in P_nA$, which are independent
modulo $P_{n+1}A$. In particular,
$\langle x \rangle \cap P_{n+1}A = 1$. By the definition of $(P_nG)$ as
an $\cal M$--universal structure, it now follows that
$P_nG$ contains infinitely many elements $x$ such that
$\langle x \rangle \cap P_{n+1}G = 1$.

\begin{theorem}\label{equal}
For all $x, y \in P_nG \setminus P_{n+1}G$,
$x \in \langle y \rangle ^{P_0G}$.
\end{theorem}

\noindent
{\bf Proof.} If $x$ and $y$ satisfy the assumptions in
Theorem \ref{conjugates}, then we are done.
Suppose we have an arbitrary element $x \in P_nG \setminus P_{n+1}G$.
We show that there exists an element
$a_x \in \langle x \rangle ^{P_0G}$ such that the elements
$a_x, xa_x \in P_nG \setminus P_{n+1}G$ have infinite order, and
$(\langle a_x \rangle \cup \langle xa_x \rangle ) \cap P_{n+1}G = 1$.
For then, by Theorem \ref{conjugates},
$$\langle xa_x \rangle ^{P_0G} = \langle a_x \rangle ^{P_0G} =
\langle a_y \rangle ^{P_0G} = \langle ya_y \rangle ^{P_0G}\,,$$ and so
$x \in \langle a_x \rangle ^{P_0G} = \langle ya_y \rangle ^{P_0G} \subseteq
\langle y \rangle ^{P_0G}$.

Suppose for now that the order of $x$ is infinite. There exists a
free group $F$ of infinite rank and contained in $P_0G$ such that
$\langle F, x \rangle = F \ast \langle x \rangle$. Let
$P_0H = F \ast \langle x \rangle$ and $P_nH = (F^{(n)})^{\langle x \rangle}$,
where $F^{(n)}$ denotes the $n$th derived group of $F$. By
Lemma \ref{intersection}, $(F^{(n)}) \leq (P_nH)$. Since $(P_nG)$ is
$\cal M$--universal, $(P_nH)$ embeds in $(P_nG)$. Without loss of generality,
assume $(P_nH) \leq (P_nG)$. Let $a \in F^{(n)} \setminus F^{(n+1)}$.
Then $x^ax \in P_nG$ of infinite order.
If $(x^ax)^k \in P_{n+1}G$, for some $k > 0$, then
$x^ax = x^2(a^{-x^2}a^x)$ and
$$(x^ax)^k = (a^{-x^2}a^x)^{x^{-2}} \cdots (a^{-x^2}a^x)^{x^{-2k}}x^{2k}
\in P_{n+1}G \cap P_0H = P_{n+1}H \,.$$
If $n \geq 1$, then $x^{2k} \in P_nH$, since $P_{n+1}H \subseteq P_nH$ and
$a \in P_nH$. By the definition of $P_nH$, $x^{2k} = 1$.
This gives a contradiction, since $x$ was assumed to have infinite order.
If $n = 0$, $(x^ax)^k \in P_1H \leq (P_0H)^{(1)}$. Since
$P_0H/(P_0H)^{(1)}$ is free abelian (see \cite{LS}, p. 24),
$$(x^ax)^k \equiv x^{2k} \equiv 1 \mbox{ mod } (P_0H)^{(1)} \,.$$
This gives a contradiction, since $x(P_0H)^{(1)}$ has infinite order as
an element of the quotient group $P_0H/(P_0H)^{(1)}$.
Similarly, the element $x\cdot x^ax \in P_nG$ has infinite order and
$\langle x\cdot x^ax \rangle \cap P_{n+1}G = 1$. Thus we can take
$a_x = x^ax$.

If the order of $x$ is finite, choose an element $g \in P_0G$ such that
$\langle x, g \rangle = \langle x \rangle \ast \langle g \rangle$.
Then $x^gx \in P_nG$ has infinite order. If
$x^gx \not\in P_{n+1}G$, apply the preceding to $x^gx$.
If $x^gx \in P_{n+1}G$, then $xx^gx \in P_nG \setminus P_{n+1}G$ has infinite
order and apply the preceding to $xx^gx$.
$\hfill \square$

\begin{theorem}
Let $I = \cap _{n < \omega}P_nG$. If $H$ is a normal subgroup of $P_0G$
such that  $H$ is not contained in $I$, then $H = P_nG$ for some $n < \omega$.
Thus every non-trivial normal subgroup of $P_0G/I$ is isomorphic to
$P_kG/I$ for some $k < \omega$.
\end{theorem}

\noindent
{\bf Proof.} Let $H$ be a normal subgroup of $P_0G$ such that
$H$ is not contained in $I$. Then
there exists a least $n$ such that $H \not\subseteq P_nG$. So
$H \subseteq P_{n-1}G$ and there exists $x \in H \setminus P_nG$. By
Theorem \ref{equal}, it follows that $H = P_{n-1}G$. $\hfill \square$

\bigskip
The proof of the main theorem will now follow, when we take $G = P_0G/I$ and
$G_n = P_nG/I$ ($0 < n < \omega$).

\bigskip
\noindent
R\"udiger G\"obel \\ Fachbereich 6, Mathematik und
Informatik
\\ Universit\"at Essen, 45117 Essen, Germany \\ {\small e--mail:
R.Goebel@Uni-Essen.De}\\
Agnes T. Paras \\ Department of Mathematics \\
University of the Philippines at Diliman \\
Quezon City, Philippines \\
{\small e--mail: agnes@math01.cs.upd.edu.ph}\\
and \\ Saharon Shelah \\ Department of
Mathematics\\ Hebrew University, Jerusalem, Israel
\\ and Rutgers University, Newbrunswick, NJ, U.S.A \\ {\small
e-mail: Shelah@math.huji.ae.il}

\end{document}